# Identifying Job Satisfaction Parameters among the Employees in Higher Educational Institutions: A Mathematical Model

**Mahak Bhatia**

Department of Science
St. Xavier's College
Jaipur, India
drmahakbhatia@gmail.com

**Aled Williams**

Department of Mathematics
London School of Economics and Political Science, London, UK
a.e.williams1@lse.ac.uk

**Abstract**

The study is conducted to evaluate the job satisfaction among the administrative and teaching faculties in higher educational institutions. Many researchers have conducted studies to evaluate differences in job perception between teaching and non-teaching staff. Despite this, none of the studies have explicitly focused on developing a formal mathematical approach for analysis. Thus, this paper aims to identify the job satisfaction parameters among staff in an educational institution by using Multi Criteria Decision Making (MCDM) tools.

The factors influencing employee job satisfaction have been ascertained and hierarchically organized through the utilization of the standard deviation methodology. Analytical findings reveal that certain variables, including promotional opportunities, interpersonal relations with colleagues, managerial support, and the department of employment, exert a substantial impact on an individual's level of job satisfaction. The research posits that the strategic implementation of these variables and attributes by organizational management can significantly ameliorate challenges related to employee retention, thereby enhancing overall workforce efficiency.

## 1. Introduction

The retention of employees in an organization is a challenging issue in the present dynamic working environment. After the outbreak of coronavirus pandemic, both the retention and recruitment of staff have grown particularly problematic for educational institutions. The incorporation of enhanced conceptualization and contextualization in educational methodologies introduced additional layers of complexity to the sector. The growing integration of technology in education is contributing to perceptions of job dissatisfaction among faculty. This is largely because many institutions have not implemented faculty development programs to equip educators with the skills needed for these new teaching approaches. Additionally, the financial disparity between teaching and non-teaching staff exacerbates feelings of job dissatisfaction. Bowen and Radhakrishna (1991) found that the level of job satisfaction among the employees is independent of their age, experience and type of institution. In contrast, Kalleberg and Loscocco (1983) found that age is positively related to job satisfaction.

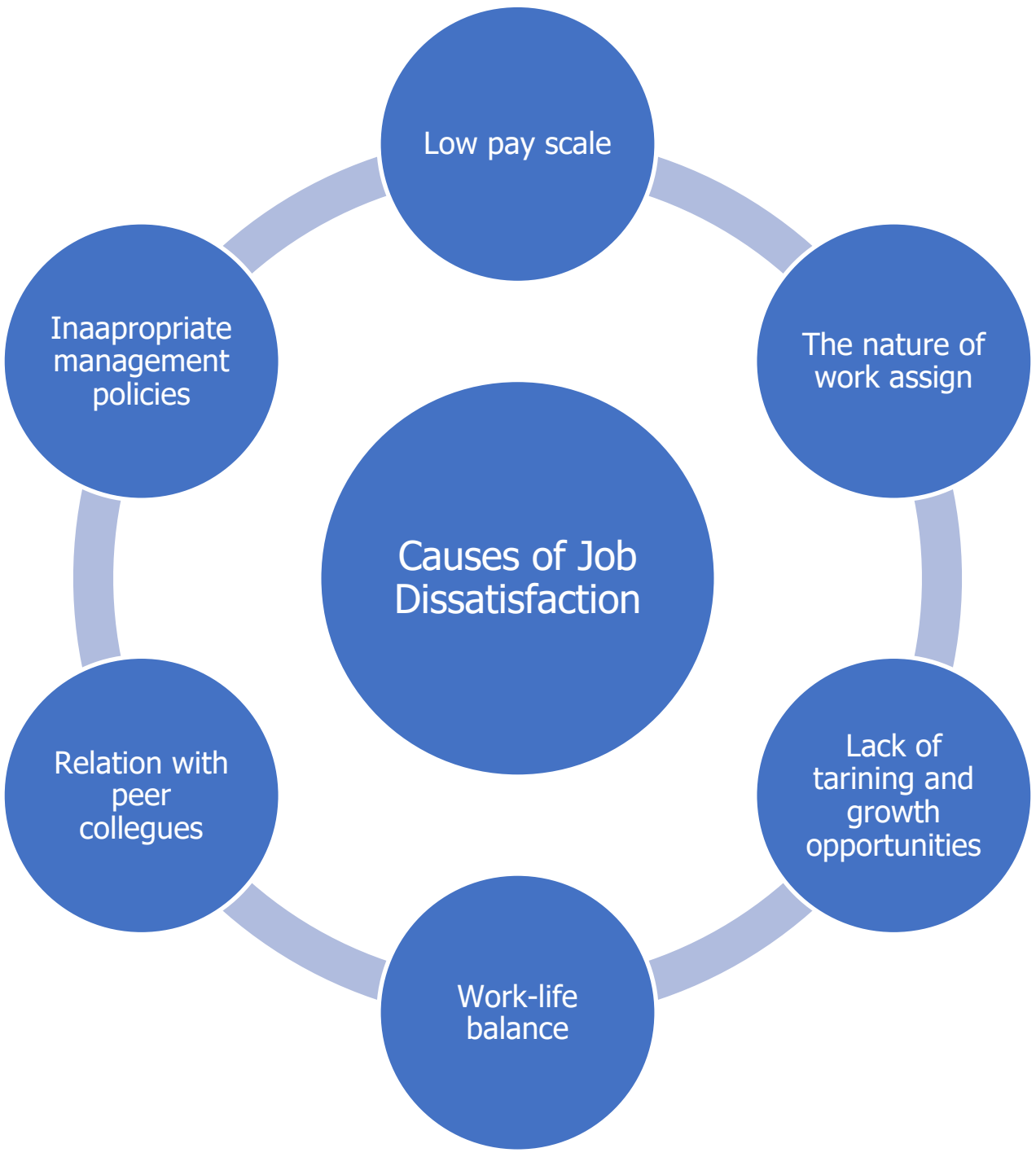

Fig 1: Parameters of job dissatisfaction

Bowen and Radhakrishna (1991) employed Herzberg's theory to assess the role of motivational factors and interpersonal relationships among colleagues in determining levels of job satisfaction. The motivational factor reveals that academics frequently derive greater satisfaction from the intrinsic qualities of their work rather than from the opportunities provided to them. Regarding relationships with peer colleagues, academics are predominantly motivated by the conduct of their co-workers and are least content with the criteria employed to determine their compensation structure.

The conclusions drawn from prior research suggest that disparities in job satisfaction may exist based on the employee's level and the nature of their assigned tasks. Locke (1976) defines job satisfaction as a pleasurable emotional state that is a result of mutual compatibility among the employees and the perception of being satisfied with the job. Job satisfaction is defined as one of the aspects that reflects employers' attitude towards their job and its related attributes. Othman (2014) found that with increases in occupational level, the job satisfaction also increases. The higher the academic position, the more satisfaction is found generally among employees. It can be therefore concluded that the role one holds within an organization also serves as a critical factor in determining job satisfaction.

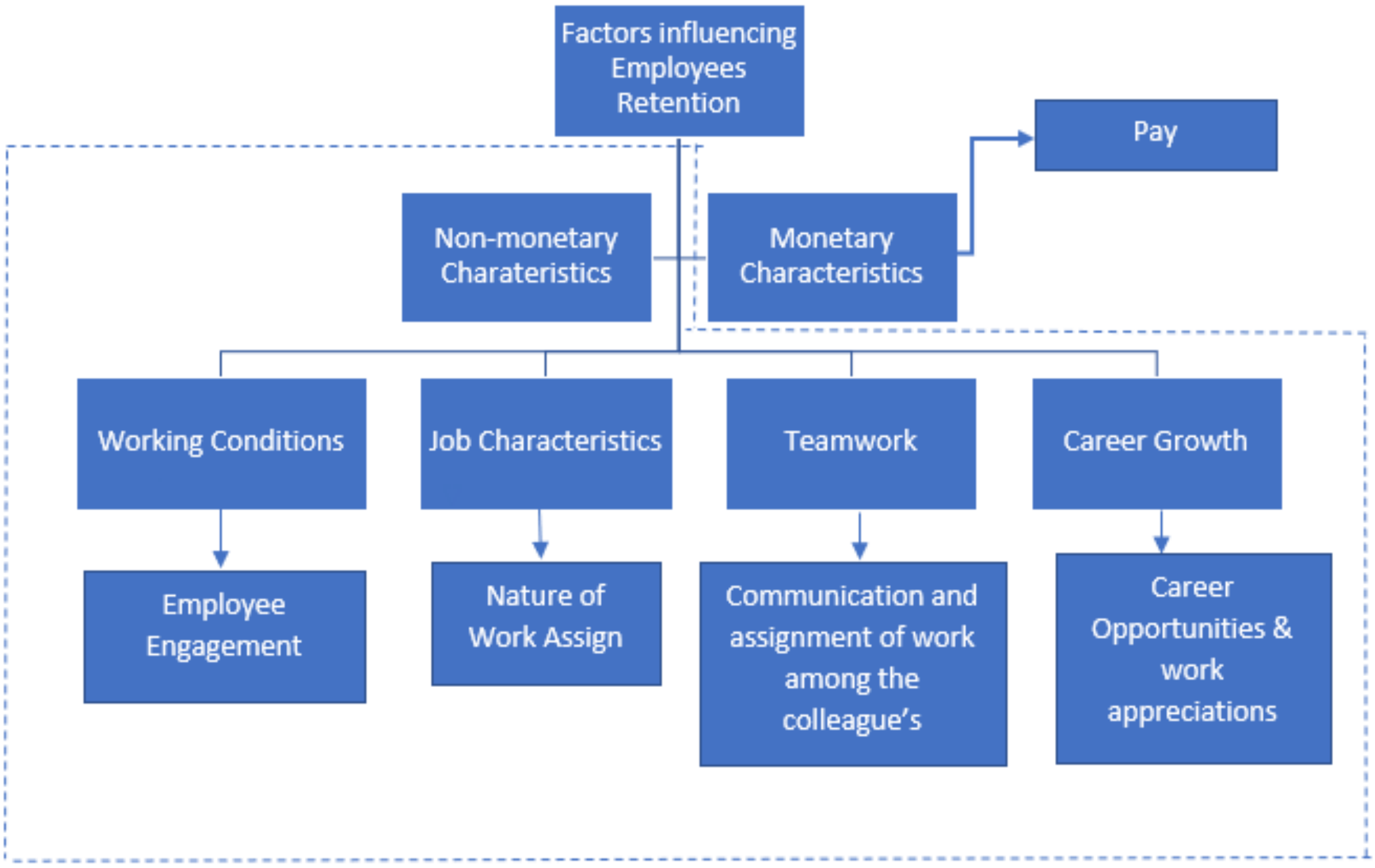

Fig.2: Monetary and Non-monetary parameters influencing employee's retention

When employees collaborate, they frequently engage in challenges, new opportunities and innovative thinking. These interactions not only make the atmosphere positive, but also this contributes to their overall job satisfaction level. Lu et al. (2016) argue that work engagement is a personal attribute, distinct to each individual. They further argue that job satisfaction is the resultant outcome derived from these individualized parameters of work engagement. In other words, the level of an employee's active involvement and enthusiasm for their work—referred to as work engagement—serves as the underlying factors that ultimately contribute to their overall job satisfaction.

It has been empirically observed that escalating work demands imposed by organizations can compromise employees' dedication to their responsibilities, subsequently diminishing their perceived levels of job satisfaction. This incongruence between individual expectations and organizational requirements often culminates in overall job dissatisfaction. Contrarily, some scholarly investigations challenge the notion that a direct correlation exists between workaholism and job satisfaction. For instance, Rayton and Yalabik (2014) conducted an in-depth exploration into the relationship among Psychological Contract Breach, work engagement, and job satisfaction. Their analytical findings indicate that employees become fully engaged in their work only when organizations fulfil all pre-established conditions. Moreover, it has been observed that when the factors affecting employee motivation are rigorously assessed, an enhancement in workforce productivity is discernible.

While ongoing research endeavours to elucidate the interrelationship between job satisfaction, work engagement, and efficiency, a more comprehensive analysis is requisite for a nuanced understanding of the variables and factors that influence job satisfaction (Garg et al., 2018). Rothmann (2008) conducts an empirical examination of the correlation between job satisfaction and work efficiency, identifying these elements as pivotal indicators of work-related well-being.

The purpose of this paper is to study the level of job satisfaction among the teaching and non-teaching faculty in higher educational organization and to analyse the relationship between the job satisfaction and work engagement. Thus, to evaluate the people's perception working in an educational institution and to identify the parameters that leads to job dissatisfaction a mathematical model is formulated that identifies the constraints that are being faced by employees. Different MCDM tools are employed that provide insight into the parameters that impact job satisfaction. To the best of our understanding, there exists no study specifically aimed at quantitatively assessing job satisfaction among employees in higher educational institutions through a mathematical methodology. The findings of this research contribute to the identification of policies and attributes that influence faculty satisfaction levels, thereby offering a scientific foundation for these determinants.

The majority of existing research posits that job satisfaction serves as a potential predictor of absenteeism, turnover, and workaholism, contending that these factors fall within the sphere of managerial influence and thus establish benchmarks for job satisfaction (Oshagbemi, 2003). Despite this, a counterargument challenges this notion by highlighting the heterogeneity in individual expectations from employment. Consequently, the mathematical model developed in this study offers nuanced insights into the diverse factors that impact employee job satisfaction.

The analysis furnishes a scientific framework for comprehending the variations in levels of job satisfaction. It is anticipated that such scholarly contributions will assist organizational management in pinpointing deficiencies and in refining policies, rules, and regulations to enhance employee motivation and engagement within institutions and organizations.

While there exists a plethora of research publications on job satisfaction, there is a notable paucity of studies that specifically focus on the job satisfaction of both teaching and non-teaching faculty in higher educational institutions through a quantitative methodology. Accordingly, the present study undertakes a comprehensive review of the literature concerning factors such as age, working hours, peer behaviour, advancement opportunities, and overall job satisfaction, aiming to identify additional constraints that influence employees' perceptions of job satisfaction.

## 2. Literature review

The literature generally indicates a strong correlation between employees possessing advanced skills or knowledge and elevated levels of job satisfaction. However, Clark (1996) and Heywood et al. (2002) assert that individuals with higher educational attainment tend to report lower levels of job satisfaction. Ward and Sloane (2000) corroborate that while earnings comparisons do play a role in shaping job satisfaction, non-pecuniary factors such as collegial relationships and the nature of the job—whether teaching or non-teaching—also serve as significant determinants. Kifle and Hailemariam (2012) examine gender-based variations in job satisfaction, considering variables like working hours, workload, and organizational support. Their findings reveal that females generally express greater satisfaction with their job roles compared to their male counterparts. Contrarily, Sloane and Ward (2001) conclude that females under the age of 35 tend to be less satisfied with their employment.

Numerous scholars, including Oshagbemi (2003), have investigated the relationship between variables such as age, gender, rank, and length of service as they relate to job satisfaction. The analysis suggests that tenure within an institution positively correlates with levels of job satisfaction. Furthermore, tenure is significantly associated with determinants of satisfaction, including performance appraisals, promotional opportunities, and a supportive learning environment. Rhodes (1983) utilized bivariate statistical analysis, revealing that age consistently exerts a positive influence on job satisfaction and a negative impact on turnover intentions. Conversely, multivariate analyses indicate that a multitude of independent variables positively correlate with both job satisfaction and job involvement, although organizational commitment shows less consistency.

Jarinto et al. (2019) applied structural equation modeling to assess job-related factors affecting employee retention in higher educational institutions. Clark (1996) delineated four principal metrics impacting job satisfaction, noting that individuals in academic and healthcare sectors are generally more satisfied by the services they provide than by their compensation. Contrary to prior research, the study found that married individuals report lower levels of satisfaction compared to their unmarried counterparts. Workplace training was found to enhance all metrics of job satisfaction.

Siregar and Siregar (2018) employed the Analytical Hierarchy Process (AHP) to rank job satisfaction parameters using a weighted mean approach. Their findings indicate that customer awareness criteria carry the highest weight, while incentives for employees exhibit variability. Mahdavi et al. (2011) introduced a fuzzy logic model to evaluate job satisfaction levels within organizations.

Through a comprehensive review of existing literature, several research gaps were identified, thereby aiding in the formulation of research questions and objectives:

1. Identification of parameters that impact the perception of job satisfaction.
2. MCDM analysis is conducted to evaluate the impact of parameters.

3. **Results**

Table 1: The alternatives and the parameters impact the job satisfaction

| **Alternatives** | **Average working hours** | **Satisfy with the pay scale with respect to the work load** | **Do you get ample opportunities at workplace to develop a skill** | **Satisfied with the working environment in an organization** | **Satisfied by the appraisals given by management** | **Satisfied with the nature of work allotted** | **Get the appreciation of the work/tasks conducted** | **Satisfied with the behaviour of peer employees in an organization** | **Satisfied with the policies and rules & regulation by the management** | **Satisfy with the designation allotted in an organization** | **Seeking to change the job if got a high pay scale** |
|---|---|---|---|---|---|---|---|---|---|---|---|
| 1 | 3.53 | 3.64 | 4.09 | 3.82 | 3.91 | 3.48 | 3.79 | 3.23 | 3.8 | 3.93 | 4.13 |
| 2 | 4.08 | 3.8 | 3.52 | 3.61 | 3.46 | 3.28 | 3.3 | 3.94 | 3.66 | 3.59 | 3.84 |
| 3 | 3.5 | 4.11 | 3.92 | 3.66 | 4.05 | 3.87 | 3.84 | 4.12 | 3.92 | 4.11 | 3.3 |
| 4 | 3.35 | 3.76 | 3.15 | 3.31 | 3.94 | 3.62 | 3.79 | 3.44 | 3.46 | 4.07 | 3.93 |
| 5 | 4.13 | 3.3 | 3.94 | 3.8 | 3.96 | 4.13 | 3.55 | 3.43 | 3.2 | 0.92 | 0.82 |
| 6 | 1.03 | 0.95 | 0.92 | 0.91 | 1.25 | 1.04 | 0.95 | 1.16 | 0.95 | 1.01 | 1.11 |
| 7 | 1 | 0.99 | 0.93 | 1.05 | 0.96 | 0.93 | 1.02 | 0.88 | 1.22 | 1.15 | 0.91 |
| 8 | 1.23 | 1.01 | 0.8 | 1.32 | 1.67 | 0.89 | 0.98 | 0.96 | 0.91 | 0.93 | 0.92 |
| 9 | 0.96 | 0.96 | 1.04 | 1.2 | 1.67 | 0.87 | 0.96 | 0.98 | 0.95 | 0.97 | 0.91 |
| Criteria | NB | NB | NB | B | B | B | B | B | B | B | B |

NB*: Non-Beneficiary

B*: Beneficiary

Table 2: Different ranks assigned to the alternatives that impact the job satisfaction level

| Alternatives | Si- | Si+ | ci | rank |
|---|---|---|---|---|
| **Average working hours** | 0.089602 | 0.055112 | 0.619168 | 3 |
| **Satisfy with the pay scale with respect to the work load** | 0.08404 | 0.054844 | 0.605111 | 4 |
| **Do you get ample opportunities at workplace to develop a skill** | 0.091001 | 0.055916 | 0.619405 | 2 |
| **Satisfied with the working environment in an organization** | 0.088131 | 0.04766 | 0.64902 | 1 |
| **Satisfied by the appraisals given by management** | 0.071592 | 0.078387 | 0.477346 | 6 |
| **Satisfied with the nature of work allotted** | 0.058143 | 0.091976 | 0.387313 | 7 |
| **Get the appreciation of the work/tasks conducted** | 0.057905 | 0.092936 | 0.383884 | 9 |
| **Satisfied with the behaviour of peer employees in an organization** | 0.058012 | 0.092506 | 0.385416 | 8 |
| **Satisfied with the policies and rules & regulation by the management** | 0.057692 | 0.092693 | 0.383628 | 10 |
| **Satisfy with the designation allotted in an organization** | 0.047449 | 0.087666 | 0.351176 | 11 |
| **Seeking to change the job if got a high pay scale** | 0.060185 | 0.065686 | 0.478147 | 5 |

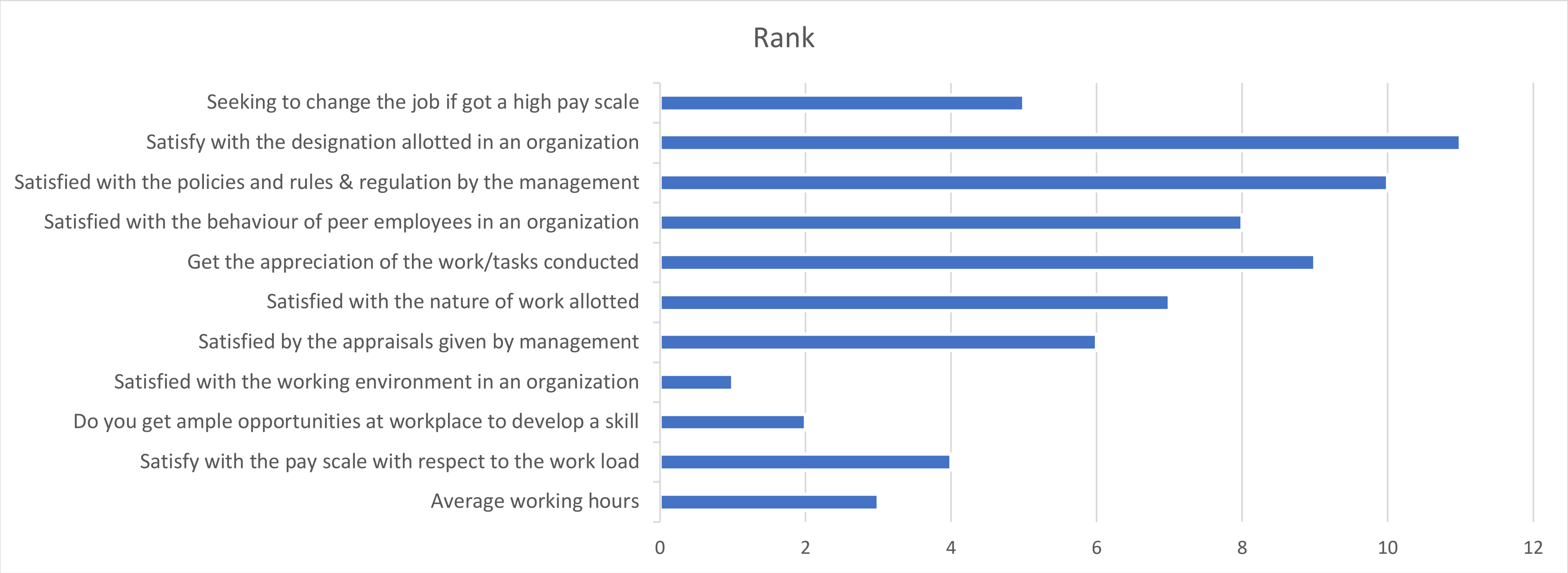

Graph: Rank assigned to different parameters of Job Satisfaction

The research offers an examination of the factors influencing job satisfaction among both teaching and non-teaching staff in higher educational institutions. Data were gathered from five distinct academic organizations utilizing questionnaires. The findings reveal that employees in these higher educational settings generally express satisfaction with their working conditions, yet display lower levels of contentment with their assigned roles and the managerial policies in place.

The survey revealed that non-teaching staff are regarded as foundational elements within academic institutions. Despite being a crucial and highly valued resource within these organizations, it was observed that the compensation structures and managerial policies are not particularly favourable to non-academic personnel. As a result, retaining such staff presents a significant challenge for educational institutions. The analysis suggests that factors like relationships with colleagues and opportunities to demonstrate skills contribute positively to their job satisfaction. However, elements such as assigned roles and compensation levels adversely affect their work efficiency.

The findings of the study have been evaluated by management professionals within academic institutions. This evaluation offers them insights into the policies that contribute to job dissatisfaction among employees. Hannay and Northam (2000) identified four key factors that significantly predict employee retention: the perception of future growth opportunities, the age of the employee, the alignment between employer and employee expectations, and the support or training provided by the employer.

The analysis indicates that there is an inverse relationship between educational qualification and levels of job satisfaction. During a survey it is analysed that that the women are more satisfied with their job profile as compared to men. The MCDM techniques used for evaluation focuses on ranking the different parameters that impact overall job satisfaction. Based on the analysis, it can be concluded that teaching staff generally exhibit higher levels of satisfaction compared to non-teaching staff, particularly in the areas of working conditions and compensation. Further, it seems like the primary cause of job dissatisfaction among the non-teaching staff is inadequacy of the salaries offered. Various MCDM techniques offer distinct perspectives on the criteria being evaluated. Therefore, consistent results for the criteria under consideration can only be expected when using a specific analytical tool. Employing different tools for analysis is likely to yield varying outcomes.

Funding: No funding is available for the research work

Disclosure Statement: There is no conflict of interests

Data availability statement; the data used for the analysis is collected from primary source